\documentclass[3p,times]{elsarticle}
\usepackage[T1]{fontenc}
\usepackage{amsmath}
\usepackage[colorlinks,linkcolor=blue,anchorcolor=red,citecolor=green]{hyperref}
\allowdisplaybreaks[4]
\usepackage{xcolor}
\usepackage{booktabs}
\usepackage{multirow} 
\usepackage{graphicx}
\usepackage{float}

\usepackage{amssymb}
\usepackage{amsthm}
\usepackage{subfigure}
\usepackage{longtable}
\usepackage{array}
\usepackage{pdflscape}
\newtheorem{theorem}{Theorem}%[section]
\newtheorem{lemma}{Lemma}%[section]
\newtheorem{proposition}{Proposition}%[section]

\usepackage[figuresright]{rotating}
\begin{document}
\begin{frontmatter}
	\title{Exact maximal order of a weak focus and a codimension-four Hopf bifurcation in a Holling--Tanner model with harvesting}%\tnoteref{label0}} %A codimension-four Hopf bifurcation and four coexisting limit cycles in a Holling--Tanner model with harvesting}
	
	%\tnotetext[label0]{This work is partially supported by the National Natural Science Foundation of China (Nos. 12071162 and 11701191), the Natural Science Foundation of Fujian Province (No. 2025J01168), Quanzhou Municipal Special Program for Key Sci-tech Talents (No. 2026QZC29GG), and the Fundamental Research Funds for the Central Universities (No. 26JCPY004).}

	\author[label1]{Zhenshu Wen}
	\ead{wenzhenshu@hqu.edu.cn}
	\author[label1]{Qiyi Wu}
	\ead{25011070006@stu.hqu.edu.cn}
	\author[label2]{Yonghui Xia\corref{cor1}}
	\ead{xiadoc@163.com}
	\cortext[cor1]{Corresponding author.}
	\address[label1]{School of Mathematical Sciences, Huaqiao University, Quanzhou 362021, PR China}
	\address[label2]{School of Mathematics, Foshan University, Foshan 528000, PR China}
	\begin{abstract}
          %Exact real-root isolation, polynomial remainder sequences, and reduction modulo a univariate elimination polynomial certify the focal-value identities, admissibility, and full-rank transversality in the original four-parameter family. Hence the system admits a codimension-four Hopf bifurcation, whose suitable generic unfolding can produce four small-amplitude limit cycles. Numerical continuation detects three folds of periodic orbits. A separate high-precision computation indicates four nested cycles with alternating stability, including two attracting cycles. The latter configuration is numerical evidence, not an interval-arithmetic proof.
          We study degenerate Hopf bifurcations in a Holling--Tanner model. Previous work established that the maximal weak-focus order lies between three and four and that a codimension-three Hopf bifurcation occurs, leaving the attainability of order four open. We prove that the maximal order is exactly four by constructing an admissible point at which the first three focal values vanish and the fourth is nonzero. Exact real-root isolation, polynomial remainder sequences, and reduction modulo a univariate elimination polynomial certify the focal-value identities, admissibility, and full-rank transversality in the original four-parameter family. Hence the system admits a codimension-four Hopf bifurcation, whose suitable generic unfolding can produce four small-amplitude limit cycles. Numerical continuation further detects three saddle-node bifurcations of limit cycles as the harvesting rate varies. We also present a parameter set for which four nested limit cycles coexist and determine their alternating stability by means of the Poincar\'e return map and the nontrivial Floquet multipliers. In this configuration, the four cycles alternate in stability, with two attracting and two unstable cycles.
    \end{abstract}	
\end{frontmatter}

\section{Introduction}\label{Intro}

A central problem in the local bifurcation theory of planar systems is to determine the order of a weak focus and its cyclicity under perturbations. At a Hopf equilibrium, the weak-focus order is determined by the first nonvanishing focal value of the local Poincar\'e return map, whereas the number of small-amplitude limit cycles realized in a parameter family also depends on the available unfolding directions. Bautin's classical
work \cite{bautin1952} established the sharp bound of three small-amplitude limit cycles bifurcating from a quadratic focus or center. Subsequently, Zoladek \cite{zoladek1994quadratic,zoladek1994generalization} revisited the quadratic case and obtained related cyclicity results for certain cubic systems. Qiu and Yang \cite{qiu2009focus} investigated the realization of high-order weak foci in polynomial vector fields, while Garc\'{\i}a, Llibre, and Maza \cite{garcia2013cyclicity} studied the cyclicity of analytic foci by means of the vanishing multiplicity of inverse integrating factors. The algebraic structure of local cyclicity was further explored through reduced
Bautin indices by Hauser, Risler, and Teissier \cite{hauser1999reduced}, and through computational algebra applied to nonradical Bautin ideals by Levandovskyy, Romanovski, and Shafer \cite{levandovskyy2009nonradical}. Related, but distinct, limit-cycle counting problems arise from perturbations of period annuli: Gavrilov, Gin\'e, and Grau \cite{gavrilov2009cyclicity} derived bounds for weight-homogeneous centers, and Gin\'e, Grau, and Llibre \cite{gine2015quasihomogeneous} established upper bounds for cycles bifurcating from quasi-homogeneous centers. In another setting, Gin\'e, Grau, and Llibre
\cite{gine2011polynomial} proved a sharp upper bound for the number of polynomial limit cycles of a class of scalar polynomial differential equations. Although these results concern different systems and notions of cyclicity, they emphasize a common mathematical challenge: obtaining sharp bounds and determining when those bounds are attained. Many notable contributions in this field concern sharp bounds for limit cycles related to the 16th Hilbert problem; see, for instance, \cite{chenLlibreWangZhang2017cyclicity,liuXiao2020abelian,chenFengJiaMeng2026maximum,chenFangZhang2025algebraic,hanZhang2010nonsmooth,chenLiYu2024nilpotent,chenTangXiao2023uniqueness,chenLlibreTang2021higgins,ruanTangZhang2010versal,xiongChen2025homoclinic}.

For rational planar systems with multiple parameters, focal values can be high-degree rational functions. Thus, rigorous determination of the exact weak-focus order and its unfolding remains a challenging problem in the analysis of nonlinear population models. This task requires careful verification of focal-value equations, nonvanishing conditions, parameter admissibility, and transversality.
In fact, several interesting results concerning limit cycles have been reported in earlier literature on population systems with harvesting strategies and specific parameter regimes; see, for example, \cite{may1979management,beddington1980maximum,brauer1979stability,brauer1982coexistence,brauer2012mathematical,dai1998coexistence,xiao1999bogdanov,xiao2005bifurcations,xiang2022degenerate,pan2024effects,gasull1997limit,xiang2023refuge,wu2026,Jia2026}.
The present paper investigates how constant-yield prey harvesting affects high-order Hopf degeneracies within this model, giving rise to remarkably rich dynamical behavior. Our motivation originates from the work of Xiang, Lu, and Huang \cite{xiang2022degenerate} on the constant-yield harvested Holling-Tanner model. 
  They established that the maximal order of a positive weak focus is at least three and at most four, and showed that the model can undergo a codimension-three Hopf bifurcation. They also obtained a codimension-four Bogdanov--Takens bifurcation. The attainability of the upper bound four for the weak-focus order, however, remained unresolved. It is important to distinguish the two codimension-four phenomena: a Hopf degeneracy is associated with a pair of nonzero purely imaginary eigenvalues, whereas a Bogdanov--Takens singularity has a double-zero eigenvalue. Consequently, the codimension-four Bogdanov--Takens result in \cite{xiang2022degenerate} does not settle whether the same model admits a fourth-order weak focus. Nor does the existence of a codimension-three Hopf point imply that the next focal degeneracy is attainable in the original four-parameter family. To settle the open question, one must construct an admissible positive Hopf equilibrium at which the first three focal values vanish, the fourth does not vanish, and the corresponding three focal-value conditions satisfy a full-rank transversality condition on the Hopf surface.

The principal result of this paper is that the upper bound is attained: the maximal order of a positive weak focus in the harvested Holling--Tanner model is exactly four. We construct an admissible point on the Hopf surface satisfying
\[
\sigma_1=\sigma_2=\sigma_3=0,\qquad \sigma_4\ne0,
\]
and prove that the Jacobian of the first three focal values with respect to coordinates on the Hopf surface is nonsingular at this point. We also show that the trace-zero condition is transverse to the remaining parameter direction and that the change from the equilibrium coordinate to the original harvesting-model parameter is locally invertible. These facts establish the codimension-four Hopf nondegeneracy in the original four-parameter family, rather than only in a reduced or artificially enlarged family. In particular, the previously known bound \(3\leq m_{\max}\leq4\) is sharpened to \(m_{\max}=4\).

The proof relies on exact algebraic certification. The elimination framework in \cite{xiang2022degenerate} already uses resultants and pseudo-division to analyze the algebraic varieties defined by the focal values. Our contribution is to complete the construction of an admissible fourth-order Hopf point through a sequence of exact verification steps. The relevant multivariate conditions are reduced to a univariate elimination polynomial. A real zero is isolated in a rational interval by exact root-counting, and polynomial remainder sequences are then used to reconstruct the remaining equilibrium and model parameters. Reduction modulo the elimination polynomial verifies the original focal-value identities, while resultant and greatest-common-divisor calculations certify the nonvanishing of the denominators, leading coefficients, and transversality determinant used in the reconstruction. Rational interval estimates establish that the reconstructed point satisfies the strict inequalities defining the admissible parameter domain. This procedure is potentially applicable to other polynomial or rational planar systems when their focal-value conditions admit tractable algebraic elimination and when the relevant inequalities can be certified exactly. Its scope is necessarily limited by the complexity of the elimination polynomials and does not extend automatically to arbitrary smooth or nonalgebraic systems.

The fourth-order weak focus has a direct local dynamical implication. The verified nondegeneracy and transversality conditions imply that a suitable generic unfolding of the original four-parameter family can produce four small-amplitude limit cycles near the positive equilibrium. We complement this local result with a numerical investigation of periodic orbits. Numerical continuation with respect to the harvesting rate detects three saddle-node bifurcations of limit cycles. At a selected parameter set, the Poincar\'e return map has four computed positive fixed points, corresponding to four nested periodic orbits. The associated nontrivial Floquet multipliers give the stability sequence, from the innermost to the outermost cycle, as unstable, stable, unstable, and stable. Thus, two attracting periodic orbits coexist, with unstable cycles separating their basins in the nested configuration. This example illustrates a form of oscillatory multistability in a planar predator--prey system. We stress that the codimension-four Hopf point is established by the exact algebraic argument, whereas the four-cycle configuration and the three cycle folds are supported by high-precision numerical computations. %and are not claimed as interval-arithmetic existence proofs.

The contributions are threefold:
\begin{enumerate}
	\item We resolve the open weak-focus-order question in \cite{xiang2022degenerate} by proving that the maximal order of a positive weak focus is exactly four and that the original four-parameter family admits a codimension-four Hopf bifurcation. %This result is distinct from the codimension-four Bogdanov--Takens bifurcation established in \cite{xiang2022degenerate}.
	\item We provide an exact algebraic certification of the fourth-order degeneracy, combining univariate elimination, rational real-root isolation, polynomial remainder sequences, reduction modulo an elimination polynomial, admissibility checks, and a full-rank transversality test. The method clarifies which algebraic nonvanishing conditions are needed to turn a candidate root into a rigorous bifurcation point.
	\item We establish the local implication of the fourth-order Hopf point for the generation of up to four small-amplitude limit cycles and numerically identify three cycle folds and a four-cycle configuration with alternating stability and two attracting periodic orbits.
\end{enumerate}

The paper is organized as follows. Section~\ref{sec:prelim} introduces the model, parameterizes the Hopf surface, and recalls the algebraic decomposition of the focal-value variety. Section~\ref{sec:certificate} constructs an admissible common zero of the first three focal-value polynomials and proves that the fourth does not vanish there. Section~\ref{sec:transversality} establishes the codimension-four transversality conditions and states the main theorem. Section~\ref{sec:numerics} presents the numerical continuation results and the four-cycle example. Section~\ref{sec:conclusion} concludes the paper.

\section{Preliminaries and algebraic decomposition}\label{sec:prelim}

This section recalls the model, the parameterization of the Hopf surface, and the algebraic decomposition of the focal-value equations obtained in \cite{xiang2022degenerate}. These results form the basis for the exact construction in Section~\ref{sec:certificate}.

To state conveniently, we use the following notation. The set $V(\xi_1,\ldots,\xi_n)$ denotes the common zero set of the polynomials $\xi_1,\ldots,\xi_n$; $\mathrm{res}(f,g,x)$ denotes the Sylvester resultant of $f$ and $g$ with respect to $x$; $\mathrm{lcoeff}(f,x)$ denotes the leading coefficient of $f$ as a polynomial in $x$; and $\mathrm{PolRem}(f,g,x)$ denotes the remainder obtained when $f$ is divided by $g$ with respect to $x$.

\subsection{The model and the Hopf surface}\label{subsec:model}

We consider the nondimensional Holling--Tanner predator--prey system with constant-yield prey harvesting:
\begin{equation}\label{eq:model}
	\begin{aligned}
		\dot x&=x(1-x)-\frac{xy}{a+x}-h,\\
		\dot y&=y\left(\delta-\frac{\beta y}{x}\right),
	\end{aligned}
\end{equation}
where $a,h,\delta$, and $\beta$ are positive parameters, and $x(t),y(t)>0$ denote the prey and predator densities, respectively. Let $E_D=(z,\delta z/\beta)$ be a positive equilibrium. Imposing the equilibrium conditions together with the trace-zero condition yields the following parameterization of the Hopf surface:
\begin{equation}\label{eq:hopf-param}
	\begin{aligned}
		a&=a_H(h,\delta,z):=\frac{z^2(1-\delta-2z)}{z(z+\delta)-h},\\[1mm]
		\beta&=\beta_H(h,\delta,z):=\frac{z\delta\,[z(z+\delta)-h]}{[z(1-z)-h]^2}.
	\end{aligned}
\end{equation}
The admissible parameter domain is
\begin{equation}\label{eq:admissible-domain}
	\Omega=\left\{(h,\delta,z)\in\mathbb{R}_+^3:0<h<\min\{z(z+\delta),z(1-z-\delta)\},\ 0<z<\frac{1-\delta}{2},\ 0<\delta<1\right\}.
\end{equation}
For every $(h,\delta,z)\in\Omega$, we have $a_H>0$, $\beta_H>0$, and $\det J(E_D)>0$. Thus, the equilibrium $E_D$ has a pair of purely imaginary eigenvalues on the trace-zero surface.

After the positive time rescaling and linear normalization introduced in \cite[Section~4]{xiang2022degenerate}, the first four focal values can be written as
\begin{equation}\label{eq:focal-values}
	\sigma_i(h,\delta,z)=K_i(h,\delta,z)f_i(h,\delta,z),\qquad i=1,2,3,4,
\end{equation}
where $K_i(h,\delta,z)\ne0$ throughout $\Omega$. Consequently,
\[
V(\sigma_1,\ldots,\sigma_m)\cap\Omega
=V(f_1,\ldots,f_m)\cap\Omega,\qquad m=3,4.
\]

For $(h,\delta,z)\in\Omega$, we say that the corresponding positive Hopf equilibrium is a weak focus of order $m$ if
\[
\sigma_1=\cdots=\sigma_{m-1}=0,\qquad \sigma_m\ne0.
\]
We denote this order by $m(h,\delta,z)$ and set
\[
m_{\max}:=\sup_{(h,\delta,z)\in\Omega}m(h,\delta,z).
\]
Since each $K_i$ is nonzero on $\Omega$, the vanishing conditions may equivalently be formulated in terms of the polynomials $f_i$.

\subsection{Algebraic decomposition of the focal-value variety}

Following the elimination procedure in \cite{xiang2022degenerate}, we eliminate the variables in the order $\delta\prec z\prec h$. The first elimination steps yield
\begin{equation}\label{eq:resultants}
	\begin{aligned}
	r_{12}&:=\operatorname{res}(f_1,f_2,h)=64z^{24}\delta^3(1+\delta-2z)(2-\delta-4z)^2(1-\delta-2z)^{10}R_1^2R_2^2g_1,\\
	r_{13}&:=\operatorname{res}(f_1,f_3,h)=-64z^{39}\delta^3(1+\delta-2z)(2-\delta-4z)^3(1-\delta-2z)^{14}R_1^2R_2^3g_2,\\
	r_{14}&:=\operatorname{res}(f_1,f_4,h)=-64z^{54}\delta^3(1+\delta-2z)(2-\delta-4z)^3(1-\delta-2z)^{18}R_1^2R_2^4g_3,\\
	r_{23}&:=\operatorname{res}(g_1,g_2,z)=c_1\delta^{73}(3+\delta)(6-\delta)^2(1-\delta)^{43}S_1S_2R_3R_4,\\
	r_{24}&:=\operatorname{res}(g_1,g_3,z)=c_2\delta^{110}(3+\delta)(6-\delta)^3(1-\delta)^{61}S_1S_2R_5R_6,
\end{aligned}
\end{equation}
where $c_1,c_2$ are nonzero rational constants. The polynomials $R_1,R_2,g_1,g_2$, and $g_3$ depend on $(\delta,z)$, whereas $S_1,S_2,R_3,R_4,R_5$, and $R_6$ depend only on $\delta$. Thus, $g_1,g_2$, and $g_3$ arise from eliminating $h$, while $R_3$ is one of the univariate factors obtained by subsequently eliminating $z$ from $g_1=g_2=0$. The expressions for $g_1$ and $R_3$, which are needed below, are given in \ref{A}; the remaining auxiliary polynomials are the same as those in \cite{xiang2022degenerate}.

The decomposition and branch exclusions established in \cite[Theorem~4.1 and Lemmas~4.2--4.3]{xiang2022degenerate} imply that
\begin{equation}\label{eq:three-focal-decomp}
	V(f_1,f_2,f_3)\cap\Omega
	=\bigl(V_{33}\cup V_{34}\bigr)\cap\Omega,
\end{equation}
where
\[
V_{33}=V(f_1,f_2,f_3,g_1,g_2,R_3),\qquad
V_{34}=V(f_1,f_2,f_3,g_1,g_2,R_4).
\]
In particular, every admissible fourth-order weak-focus point must lie on one of these two branches. The construction in this paper is carried out on the $R_3$-branch, $V_{33}\cap\Omega$.

For the simultaneous vanishing of the first four focal-value polynomials, the remaining algebraic branch is contained in
\[
V\bigl(f_1,f_2,f_3,f_4,g_1,g_2,g_3,R_3R_4,R_5R_6\bigr).
\]
The exact resultant computation
\begin{equation}\label{eq:branch-exclusion}
	\operatorname{res}(R_3R_4,R_5R_6,\delta)\ne0
\end{equation}
excludes this branch. Together with the exclusion of the other branches established in \cite{xiang2022degenerate}, this gives
\begin{equation}\label{eq:no-four-focal-zero}
	V(f_1,f_2,f_3,f_4)\cap\Omega=\varnothing.
\end{equation}

\begin{proposition}\label{prop:upper-bound}
	Every positive Hopf equilibrium of system~\eqref{eq:model} in $\Omega$ has weak-focus order at most four. In particular, $m_{\max}\leq4$.
\end{proposition}
\begin{proof}
	By \eqref{eq:no-four-focal-zero}, the first four focal-value polynomials cannot vanish simultaneously at any point of $\Omega$. Since $K_i\ne0$ on $\Omega$, the same statement holds for the focal values $\sigma_i$. Hence at least one of $\sigma_1,\ldots,\sigma_4$ is nonzero, and the first nonzero focal value has index at most four.
\end{proof}

The branch decomposition and exclusions used to obtain \eqref{eq:no-four-focal-zero} are recalled from \cite[Theorem~4.1 and Lemmas~4.2--4.3]{xiang2022degenerate}. The remaining issue is attainability: we must show that $V(f_1,f_2,f_3)\cap\Omega$ is nonempty. The next section constructs such a point exactly.

%It follows that the weak-focus order of a positive Hopf equilibrium is at most four. This upper bound alone does not establish that $V(f_1,f_2,f_3)\cap\Omega$ is nonempty. We prove this attainability in the next section.

\section{An exact algebraic certificate for a fourth-order weak focus}\label{sec:certificate}

We construct a common zero of $f_1,f_2$, and $f_3$ in $\Omega$ in three steps: first, we isolate the relevant real zero of $R_3$; next, we recover $z$ from $g_1$ and $g_2$; finally, we determine $h$ from $f_1$ and $f_2$. We use exact polynomial arithmetic throughout.

For $p\in\mathbb Q[\delta]$, define
\begin{equation}\label{eq:reduction}
	[p]_{R_3}:=\mathrm{PolRem}(p,R_3;\delta).
\end{equation}
Thus, $[p]_{R_3}$ is the unique polynomial of degree less than $\deg R_3$ representing the residue class of $p$ in $\mathbb Q[\delta]/(R_3)$. This notation also applies to a rational function $a/b$, provided $a,b\in\mathbb Q[\delta]$ and $\gcd(b,R_3)=1$. In that case, the extended Euclidean algorithm described in Section~3.2 of \cite{vonzurGathen2013modern} gives polynomials $u,v\in\mathbb Q[\delta]$ such that
\[
ub+vR_3=1,
\]
and we set
\begin{equation}\label{eq:rational-reduction}
	\left[\frac{a}{b}\right]_{R_3}:=[au]_{R_3}.
\end{equation}
In particular, if $\alpha$ is a zero of $R_3$, then $b(\alpha)\ne0$ and
\[
\left[\frac{a}{b}\right]_{R_3}(\alpha)=\frac{a(\alpha)}{b(\alpha)}.
\]
Whenever a multivariate expression appears inside the reduction brackets, all indicated substitutions are made first, so that the expression being reduced depends only on $\delta$. To control the size of intermediate expressions, we reduce modulo $R_3$ after each multiplication.

All polynomial identities, resultant calculations, and greatest-common-divisor tests below are understood over $\mathbb Q$ and are carried out by exact rational arithmetic. Decimal approximations are reported only to describe the resulting algebraic point; they are not used to establish focal-value vanishing or nonvanishing.

\begin{lemma}\label{lem:remainder}
	Let $F_0$ and $F_1$ be polynomials in $x$ whose coefficients are polynomials in the parameters. Set
	\[
	F_{j+1}=\mathrm{PolRem}(F_{j-1},F_j,x),
	\qquad
	\ell_j=\mathrm{lcoeff}(F_j,x),
	\qquad j=1,\ldots,N.
	\]
	If each $\ell_j$ is well-defined and nonzero on a subset $\mathcal B$ of the parameter-variable space, then the remainder sequence is well-defined on $\mathcal B$, and
	\[
	V(F_{j-1},F_j)\cap\mathcal B
	=
	V(F_j,F_{j+1})\cap\mathcal B,
	\qquad j=1,\ldots,N.
	\]
	Consequently,
	\[
	V(F_0,F_1)\cap\mathcal B
	=
	V(F_N,F_{N+1})\cap\mathcal B.
	\]
\end{lemma}
\begin{proof}
	At each step, polynomial division gives an identity of the form
	$F_{j-1}=Q_jF_j+F_{j+1}$ on $\mathcal B$. Hence $F_{j-1}=F_j=0$ if and only if $F_j=F_{j+1}=0$ on $\mathcal B$. Iterating this equivalence proves the result.
\end{proof}

\paragraph{Step 1. Isolating the zeros of $R_3(\delta)$}
Let $R_3(\delta)$ be the degree-34 polynomial listed in \ref{A}. Exact real-root isolation (using Maple's \texttt{realroot} command) shows that $R_3$ has exactly two zeros, denoted by $\delta_1$ and $\delta_2$, in $(0,1)$. They are isolated, respectively, by the rational intervals
\begin{equation}\label{eq:delta-intervals}
	\begin{aligned}
		I_1&=\left[\frac{1599280774005682713015}{2361183241434822606848},
		\frac{799640387002841356521}{1180591620717411303424}\right],\\[1mm]
		I_2&=\left[\frac{68666068149881493191}{73786976294838206464},
		\frac{8583258518735186649}{9223372036854775808}\right].
	\end{aligned}
\end{equation}

\paragraph{Step 2. Determining $z_*$ from the remainder sequence of $g_1$ and $g_2$}
We seek common zeros of $g_1(z,\delta)$ and $g_2(z,\delta)$ associated with the two zeros of $R_3$ in $I_1$ and $I_2$. Regard $g_1$ and $g_2$ as polynomials in $z$ of degrees 14 and 30, respectively, and set
\[
\widetilde g_0(z):=g_2(z,\delta),\qquad
\widetilde g_1(z):=g_1(z,\delta).
\]
Their polynomial remainder sequence is defined by
\[
\widetilde g_{j+1}(z):=\mathrm{PolRem}(\widetilde g_{j-1},\widetilde g_j;z),
\qquad j=1,\ldots,14.
\]
In particular,
\begin{equation}\label{eq:g15}
	\widetilde g_{15}(z):=\mathrm{PolRem}(\widetilde g_{13},\widetilde g_{14};z)
	=\frac{c_5(\delta-6)(1-\delta)^9\delta^{12}(3+\delta)S_1(\delta)S_2(\delta)R_3(\delta)R_4(\delta)
		\displaystyle\prod_{i=4}^{9}Q_i(\delta)^2}
	{\displaystyle 8\prod_{i=10}^{16}Q_i(\delta)^2},
\end{equation}
where $c_5$ is a nonzero rational constant. Here $\widetilde g_{14}$ is linear in $z$, whereas $\widetilde g_{15}$ is independent of $z$. Exact resultant computations give
\begin{equation}\label{eq:resultant-nonzero}
	\begin{gathered}
		\operatorname{res}(S_1,R_3;\delta)\ne0,\qquad
		\operatorname{res}(S_2,R_3;\delta)\ne0,\qquad
		\operatorname{res}(R_4,R_3;\delta)\ne0,\\
		\operatorname{res}(Q_i,R_3;\delta)\ne0,\qquad i=4,\ldots,16.
	\end{gathered}
\end{equation}
The leading coefficients $\mathrm{lcoeff}(\widetilde g_j,z)$, $j=1,\ldots,14$, are rational functions of $\delta$ whose numerators and denominators involve only the factors $Q_i$, $i=4,\ldots,16$, and factors that are nonzero on $V(R_3)\cap\Omega$. The resultant relations in \eqref{eq:resultant-nonzero} therefore show that these leading coefficients are well-defined and nonzero on this set. In particular,
\[
V\bigl(R_3,\mathrm{lcoeff}(\widetilde g_j,z)\bigr)\cap\Omega=\varnothing,
\qquad j=1,\ldots,14.
\]
Applying Lemma~\ref{lem:remainder} to the remainder sequence in $z$ on $\mathcal B=V(R_3)\cap\Omega$, and using the fact that $\widetilde g_{15}=0$ on this set by \eqref{eq:g15}, we obtain
\begin{equation}\label{eq:g-reduction}
	V(g_1,g_2,R_3)\cap\Omega
	=V(\widetilde g_{14},\widetilde g_{15},R_3)\cap\Omega
	=V(\widetilde g_{14},R_3)\cap\Omega.
\end{equation}
Solving $\widetilde g_{14}=0$ gives
\begin{equation}\label{eq:ztilde}
	z=\widetilde Q(\delta):=\frac{(1-\delta)^3Q_1(\delta)}{Q_2(\delta)}.
\end{equation}
Differentiating this expression yields
\[
\widetilde Q'(\delta)=-\frac{(1-\delta)^2Q_3(\delta)}{Q_2(\delta)^2}.
\]
Sturm's theorem shows that neither $Q_2$ nor $Q_3$ vanishes on $I_1\cup I_2$. Hence $\widetilde Q$ is monotone on each interval. Evaluating it at the endpoints and rounding outward gives
\begin{equation}\label{eq:z-intervals}
	\begin{aligned}
		\widetilde Q(I_1)&\subset L_1:=\left[\frac{23956749}{1000000000},
		\frac{23956749001}{1000000000000}\right],\\[1mm]
		\widetilde Q(I_2)&\subset L_2:=\left[-\frac{1258299529}{1000000000000},
		-\frac{157287441}{125000000000}\right].
	\end{aligned}
\end{equation}
Since $L_2\subset(-\infty,0)$, the zero $\delta_2\in I_2$ cannot yield an admissible point. A direct rational endpoint calculation gives
\[
2\sup L_1+\sup I_1<1.
\]
Thus every point of $I_1\times L_1$ satisfies $0<\delta<1$ and $0<z<(1-\delta)/2$. Let $\delta_*:=\delta_1$ be the unique zero of $R_3$ in $I_1$, and set $z_*:=\widetilde Q(\delta_*)$.

An exact computation gives $\gcd(Q_2,R_3)=1$. Therefore, there exist $U,V\in\mathbb Q[\delta]$ such that
\[
UQ_2+VR_3=1.
\]
Combining this identity with \eqref{eq:ztilde} and the reduction convention in \eqref{eq:rational-reduction}, define
\begin{equation}\label{eq:Z-def}
	Z(\delta):=\bigl[(1-\delta)^3Q_1(\delta)U(\delta)\bigr]_{R_3}.
\end{equation}
For every zero $\delta_0$ of $R_3$, we then have $Z(\delta_0)=\widetilde Q(\delta_0)$. In particular, $Z(\delta_*)=z_*$. Exact substitution and reduction modulo $R_3$ give
\begin{equation}\label{eq:g1g2-zero}
	\bigl[g_i(Z(\delta),\delta)\bigr]_{R_3}=0,
	\qquad i=1,2.
\end{equation}
Evaluating at $\delta=\delta_*$ yields
\begin{equation}\label{eq:g-common-zero}
	g_1(z_*,\delta_*)=g_2(z_*,\delta_*)=0.
\end{equation}

\paragraph{Step 3. Determining $h_*$ from the remainder sequence of $f_1$ and $f_2$}
Having obtained $(\delta_*,z_*)$ satisfying $g_1=g_2=0$, we now determine $h$ from $f_1=f_2=0$. Set
\[
\widetilde f_0:=f_2,\qquad \widetilde f_1:=f_1,
\]
and compute
\[
\widetilde f_{j+1}=\mathrm{PolRem}(\widetilde f_{j-1},\widetilde f_j;h),
\qquad j=1,2,
\]
together with
\begin{equation}\label{eq:f4-remainder}
	\widetilde f_4=\mathrm{PolRem}(\widetilde f_2,\widetilde f_3;h)
	=\frac{C(\delta,z)g_1(z,\delta)}{D(\delta,z)}.
\end{equation}
Here $\widetilde f_3$ is linear in $h$, while $\widetilde f_4$ is independent of $h$. Write
\[
\mathrm{lcoeff}(\widetilde f_j,h)=\frac{A_j(\delta,z)}{B_j(\delta,z)},
\qquad j=1,2,3,
\]
where $A_j,B_j,C,D\in\mathbb Q[\delta,z]$ and the fractions are in reduced form. Exact computations give
\begin{equation}\label{eq:coeff-gcd}
	\begin{aligned}
		\gcd\bigl([A_j(\delta,Z(\delta))]_{R_3},R_3\bigr)
		&=\gcd\bigl([B_j(\delta,Z(\delta))]_{R_3},R_3\bigr)=1,
		\qquad j=1,2,3,\\
		\gcd\bigl([D(\delta,Z(\delta))]_{R_3},R_3\bigr)&=1.
	\end{aligned}
\end{equation}
Since $z=Z(\delta)$ on the branch $g_1=g_2=R_3=0$ in $\Omega$, these identities imply that $A_j$, $B_j$, and $D$ are nonzero on this branch. Hence the leading coefficients $\mathrm{lcoeff}(\widetilde f_j,h)$, $j=1,2,3$, are well-defined and nonzero there. Applying Lemma~\ref{lem:remainder} to the remainder sequence in $h$ on $\mathcal B=V(g_1,g_2,R_3)\cap\Omega$, and using \eqref{eq:f4-remainder}, we obtain
\begin{equation}\label{eq:f-reduction}
	V(f_1,f_2,g_1,g_2,R_3)\cap\Omega
	=V(\widetilde f_3,g_1,g_2,R_3)\cap\Omega.
\end{equation}

Solving $\widetilde f_3=0$ gives
\begin{equation}\label{eq:htilde}
	h=\widetilde P(\delta,z):=\frac{zP_3(\delta,z)}{P_2(\delta,z)},
\end{equation}
where $P_2,P_3\in\mathbb Q[\delta,z]$. Define
\[
\widehat P_2(\delta):=[P_2(\delta,Z(\delta))]_{R_3}.
\]
Since $\gcd(\widehat P_2,R_3)=1$, there exists $W\in\mathbb Q[\delta]$ such that
\[
[W\widehat P_2]_{R_3}=1.
\]
Using $W$ to represent the inverse of $\widehat P_2$, define
\begin{equation}\label{eq:H-def}
	H(\delta):=\bigl[Z(\delta)P_3(\delta,Z(\delta))W(\delta)\bigr]_{R_3}
	=\bigl[\widetilde P(\delta,Z(\delta))\bigr]_{R_3}.
\end{equation}
It follows that $H(\delta_*)=\widetilde P(\delta_*,z_*)=h_*$. Exact substitution into the original focal-value polynomials gives
\begin{equation}\label{eq:all-three-focal}
	\bigl[f_i(H(\delta),\delta,Z(\delta))\bigr]_{R_3}=0,
	\qquad i=1,2,3.
\end{equation}
Evaluating at $\delta_*$, we obtain
\begin{equation}\label{eq:f123-zero}
	f_1(h_*,\delta_*,z_*)=f_2(h_*,\delta_*,z_*)=f_3(h_*,\delta_*,z_*)=0.
\end{equation}
In particular, the value $h_*$ recovered from $f_1=f_2=0$ also satisfies $f_3=0$.

It remains to verify that the constructed point belongs to $\Omega$. Consider the following polynomial representatives modulo $R_3$:
\[
H,\quad [Z(Z+\delta)-H]_{R_3},\quad
[Z(1-Z-\delta)-H]_{R_3},\quad
\left[\frac{1-\delta}{2}-Z\right]_{R_3}.
\]
Exact real-root counting with Maple's \texttt{realroot} command shows that none of these polynomials has a zero in $I_1$. Exact rational evaluation at the midpoint of $I_1$ shows that all four are positive there. They are therefore positive at $\delta_*$, and hence
\begin{equation}\label{eq:admissibility-check}
	\begin{aligned}
		h_*>0,\qquad
		z_*(z_*+\delta_*)-h_*>0,\qquad
		z_*(1-z_*-\delta_*)-h_*>0,\qquad
		\frac{1-\delta_*}{2}-z_*>0.
	\end{aligned}
\end{equation}
Together with $(\delta_*,z_*)\in I_1\times L_1$, these inequalities show that $(h_*,\delta_*,z_*)\in\Omega$. Numerically,
\[
\delta_*\approx0.6773217537,\qquad
z_*\approx0.02395674900,\qquad
h_*\approx0.005562357091.
\]
The two upper-bound margins in the harvesting constraint are approximately $0.01123799598$ and $0.001594038840$, respectively.

By \eqref{eq:g-common-zero} and \eqref{eq:f123-zero}, the set $V_{33}\cap\Omega$ is nonempty. On the other hand, \eqref{eq:no-four-focal-zero} rules out the possibility that $f_4$ also vanishes at this point. Thus,
\[
f_4(h_*,\delta_*,z_*)\ne0,
\]
and the corresponding positive Hopf equilibrium is a weak focus of order exactly four.

\section{Transversality of the codimension-four Hopf point}\label{sec:transversality}

We now verify that the fourth-order weak focus constructed in Section~\ref{sec:certificate} is a codimension-four Hopf point in the original parameter family.

We first use $z$ in place of $\beta$ as a local parameter near the constructed point. At the positive equilibrium $E_D=(z,\delta z/\beta)$, the equilibrium conditions give
\begin{equation}\label{eq:beta-equilibrium}
	\beta=\frac{\delta z^2}{(a+z)[z(1-z)-h]}.
\end{equation}
Evaluating the Jacobian of system~\eqref{eq:model} at $E_D$ and substituting \eqref{eq:beta-equilibrium}, we obtain
\begin{equation}\label{eq:det-trace}
	\begin{aligned}
		\det J(E_D)&=\frac{\delta\,[z^3-hz+a(z-2h)]}{z(a+z)},\\[1mm]
		tr J(E_D)&=1-2z-\frac{a[z(1-z)-h]}{z(a+z)}-\delta.
	\end{aligned}
\end{equation}
Differentiating \eqref{eq:beta-equilibrium} with respect to $z$, while keeping $a,h,\delta$ fixed, yields
\begin{equation}\label{eq:beta-derivative}
	\frac{\partial\beta}{\partial z}
	=\frac{\delta z[z^3-hz+a(z-2h)]}{(a+z)^2[z(1-z)-h]^2}
	=\frac{z^2\det J(E_D)}{(a+z)[z(1-z)-h]^2}>0
\end{equation}
at the constructed point, where $\det J(E_D)>0$. Define
\[
\Phi(a,h,\delta,z)=\bigl(a,h,\delta,\beta(a,h,\delta,z)\bigr),
\]
where $\beta(a,h,\delta,z)$ is given by \eqref{eq:beta-equilibrium}. Its Jacobian determinant is
\[
\det D\Phi=\frac{\partial\beta}{\partial z},
\]
which is positive at the constructed point by \eqref{eq:beta-derivative}. The inverse function theorem therefore shows that $(a,h,\delta,z)$ are valid local coordinates for the original parameter family $(a,h,\delta,\beta)$ near this point. 

With $h,\delta$, and $z$ fixed and $\beta$ determined by \eqref{eq:beta-equilibrium}, we have
\begin{equation}\label{eq:trace-transverse}
	\frac{\partial}{\partial a} tr J(E_D)
	=-\frac{z(1-z)-h}{(a+z)^2}<0.
\end{equation}
Thus, the Hopf surface is locally represented by $a=a_H(h,\delta,z)$, as in \eqref{eq:hopf-param}, and variation of $a$ with $h,\delta,z$ fixed is transverse to this surface.

We next verify transversality within the Hopf surface. Using the polynomials $H$ and $Z$ defined in Section~\ref{sec:certificate}, differentiate $f_1,f_2$, and $f_3$ with respect to the independent variables $h,\delta$, and $z$, and then substitute $h=H(\delta)$ and $z=Z(\delta)$. Exact computation gives
\begin{equation}\label{eq:gcd-transversality}
	\gcd\left(
	\left[\left.\det\frac{\partial(f_1,f_2,f_3)}{\partial(h,\delta,z)}\right|_{h=H(\delta),\,z=Z(\delta)}\right]_{R_3},R_3
	\right)=1.
\end{equation}
Since $R_3(\delta_*)=0$, $H(\delta_*)=h_*$, and $Z(\delta_*)=z_*$, it follows that
\begin{equation}\label{eq:focal-transversality}
	\left.\det\frac{\partial(f_1,f_2,f_3)}{\partial(h,\delta,z)}\right|_{(h,\delta,z)=(h_*,\delta_*,z_*)}\ne0.
\end{equation}
At the common zero of $f_1,f_2$, and $f_3$, the identities $\sigma_i=K_if_i$ imply $\sigma_i=K_i f_i$ on the Hopf surface for $i=1,2,3$. Since $K_i(h_*,\delta_*,z_*)\ne0$, we obtain
\begin{equation}\label{eq:sigma-transversality}
	\left.\det\frac{\partial(\sigma_1,\sigma_2,\sigma_3)}{\partial(h,\delta,z)}\right|_{(h,\delta,z)=(h_*,\delta_*,z_*)}
	=\left(\prod_{i=1}^{3}K_i(h_*,\delta_*,z_*)\right)
	\left.\det\frac{\partial(f_1,f_2,f_3)}{\partial(h,\delta,z)}\right|_{(h,\delta,z)=(h_*,\delta_*,z_*)}
	\ne0.
\end{equation}
For completeness, define the map in these local coordinates by
\[
\Psi(a,h,\delta,z)=\bigl(tr J(E_D),\sigma_1(h,\delta,z),
\sigma_2(h,\delta,z),\sigma_3(h,\delta,z)\bigr),
\]
where $\beta$ is eliminated by \eqref{eq:beta-equilibrium}. Its derivative has the block-triangular form
\[
D\Psi=
\begin{pmatrix}
	\partial_a tr J & *\\
	0_{3\times1} & D_{(h,\delta,z)}(\sigma_1,\sigma_2,\sigma_3)
\end{pmatrix}.
\]
Consequently,
\[
\det D\Psi=
\frac{\partial}{\partial a}tr J(E_D)\;
\det\frac{\partial(\sigma_1,\sigma_2,\sigma_3)}
{\partial(h,\delta,z)}\ne0.
\]
This gives the full-rank condition in the local coordinates $(a,h,\delta,z)$. The local diffeomorphism $\Phi$ transfers the same rank condition to the original parameter coordinates $(a,h,\delta,\beta)$.

Together with the transverse trace condition \eqref{eq:trace-transverse} and the local invertibility of the parameter change, this verifies the codimension-four nondegeneracy condition in the original four-parameter family.

\begin{theorem}\label{thm:main}
	For system~\eqref{eq:model}, the maximal order of a positive weak focus is exactly four. More precisely, let $\delta_*$ be the unique zero of $R_3$ in $I_1$, and set $z_*=Z(\delta_*)$ and $h_*=H(\delta_*)$. The corresponding original parameters are $a_*=a_H(h_*,\delta_*,z_*)$ and $\beta_* =\beta_H(h_*,\delta_*,z_*)$. At this admissible point,
	\[
	\sigma_1=\sigma_2=\sigma_3=0,\qquad \sigma_4\ne0,
	\]
	and the full-rank condition $\det D\Psi\ne0$ holds in local coordinates for the original four-parameter family. Consequently, a suitable generic unfolding can produce four small-amplitude limit cycles near the corresponding positive equilibrium.
\end{theorem}
\begin{proof}
	The upper bound follows from \eqref{eq:no-four-focal-zero} and the nonvanishing of the factors $K_i$ in \eqref{eq:focal-values}. The construction in Section~\ref{sec:certificate} gives an admissible point at which $f_1=f_2=f_3=0$ and $f_4\ne0$, and hence the first three focal values vanish while the fourth is nonzero. Finally, \eqref{eq:trace-transverse} and \eqref{eq:sigma-transversality} establish the required transversality in the original four-parameter family. Therefore, the maximal weak-focus order is exactly four, and the associated Hopf bifurcation has codimension four. The standard local unfolding of a fourth-order weak focus then permits four small-amplitude limit cycles.
\end{proof}

\section{Numerical investigation of periodic orbits}\label{sec:numerics}

In this section, we investigate numerically the periodic-orbit structure associated with the codimension-four Hopf point established in Sections~\ref{sec:certificate} and~\ref{sec:transversality}. We focus on the detection of saddle-node bifurcations of limit cycles and on a parameter regime in which four limit cycles coexist.

\subsection{Continuation with respect to the harvesting rate}\label{subsec:continuation}

To locate four coexisting limit cycles, we first search numerically for parameter values near the codimension-four Hopf point constructed above. We adjust the parameters until the Poincar\'e return map has four distinct positive fixed points, and then continue the corresponding periodic orbits. This numerical search leads to the following parameter values:
\begin{equation}\label{eq:numerical-parameters}
	\begin{aligned}
		a_0&=0.012815832749930476, & \delta_0&=0.6822604127797279,\\
		h_0&=0.004885786597099207, & z_0&=0.023683375784836647.
	\end{aligned}
\end{equation}
For reproducibility, these finite decimal values are taken as the defining input data in the independent high-precision computation. Here, $z_0$ is the prey coordinate of the positive equilibrium. Set
\begin{equation}\label{eq:y0-beta0}
	y_0=(a_0+z_0)\left(1-z_0-\frac{h_0}{z_0}\right),
	\qquad
	\beta_0=\frac{\delta_0z_0}{y_0}.
\end{equation}
Then $E_0=(z_0,y_0)$ is an equilibrium of \eqref{eq:model} at $(a,h,\delta,\beta)=(a_0,h_0,\delta_0,\beta_0)$, where
\[
y_0\approx0.02810514194274154,
\qquad
\beta_0\approx0.5749207661679539.
\]

\begin{figure}[!t]
\centering
\includegraphics[width=0.92\textwidth]{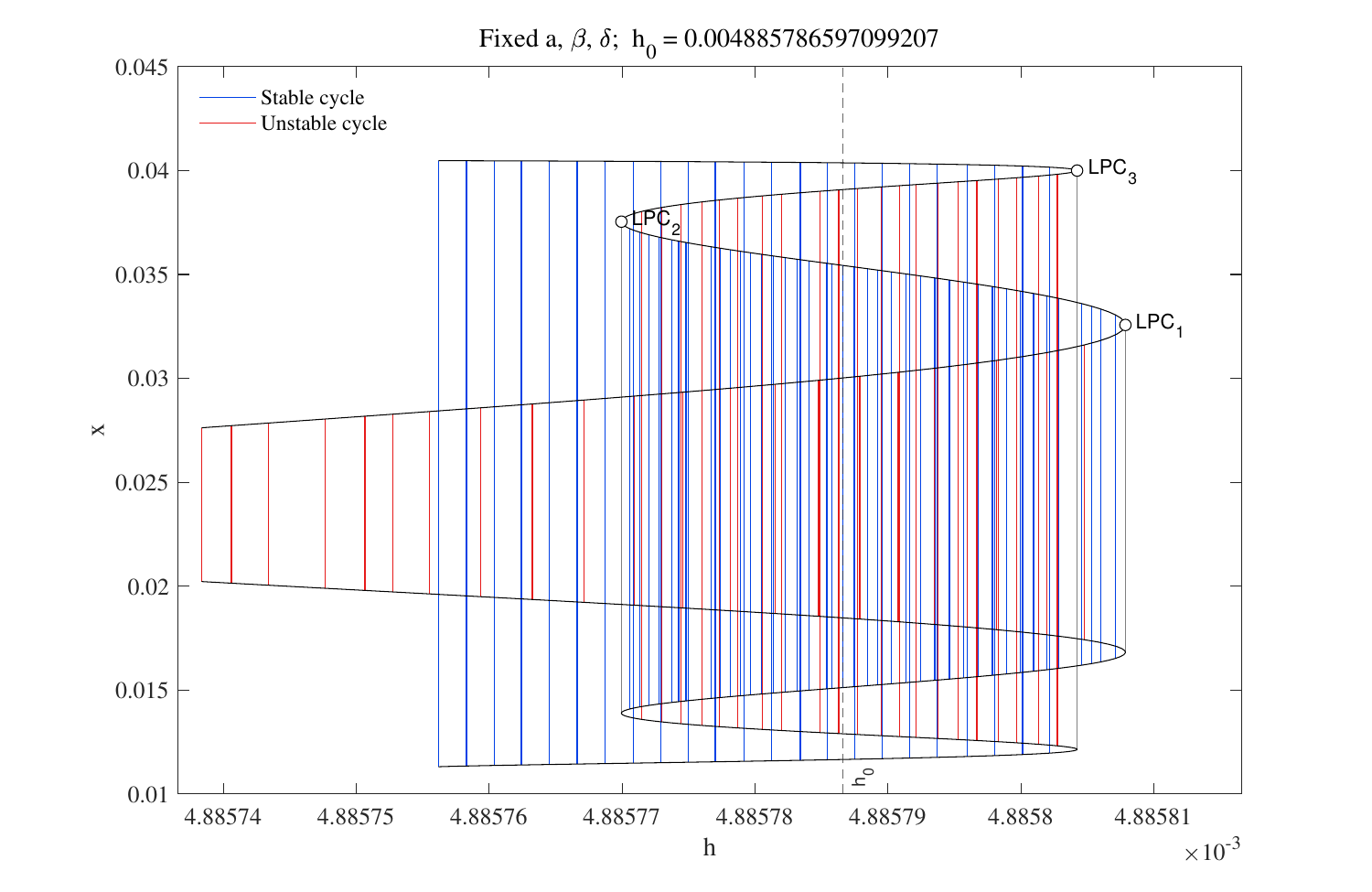}
\caption{Numerical continuation of periodic orbits of system~\eqref{eq:model} with $a=a_0$, $\delta=\delta_0$, and $\beta=\beta_0$. The black curves show the extrema of $x$, while the vertical blue and red segments represent stable and unstable cycles, respectively. The points $\mathrm{LPC}_j$ mark folds of periodic orbits. The dashed line $h=h_0$ lies in the parameter interval in which four computed cycles coexist.}
\label{fig:cycle-continuation}
\end{figure}

Fixing $a=a_0$, $\delta=\delta_0$, and $\beta=\beta_0$, we continue the periodic orbits with respect to the harvesting rate $h$ using MatCont \cite{dhooge2003matcont}. Figure~\ref{fig:cycle-continuation} shows the minimum and maximum prey densities along the computed periodic orbits. Each vertical segment represents the range of $x$ along one periodic orbit; blue and red denote stable and unstable cycles, respectively. Three saddle-node bifurcations of limit cycles, labeled $\mathrm{LPC}_1$, $\mathrm{LPC}_2$, and $\mathrm{LPC}_3$ in MatCont, are detected at approximately
\begin{equation}\label{eq:lpc-values}
	\begin{aligned}
		h_{\mathrm{LPC}_1}=0.0048858078670,\qquad
		h_{\mathrm{LPC}_2}=0.0048857699446,\qquad
		h_{\mathrm{LPC}_3}=0.0048858042316.
	\end{aligned}
\end{equation}
These approximate locations are cross-checked using the Poincar\'e return map $P_h$ by solving the fold equations
\begin{equation}\label{eq:lpc-return}
	P_h(r)-r=0,\qquad \frac{\partial P_h}{\partial r}(r)-1=0,
\end{equation}
where $r$ is a coordinate on a transverse section. The continuation calculation also indicates the nondegeneracy of the detected folds; all these fold statements are numerical. The computed bifurcation diagram exhibits four coexisting limit cycles for
\[
h_{\mathrm{LPC}_2}<h<h_{\mathrm{LPC}_3},
\]
and the selected value $h_0$ lies in this interval.

\subsection{Four coexisting cycles and their stability}\label{subsec:four-cycles}

\begin{figure}[!t]
\centering
\includegraphics[width=0.49\textwidth]{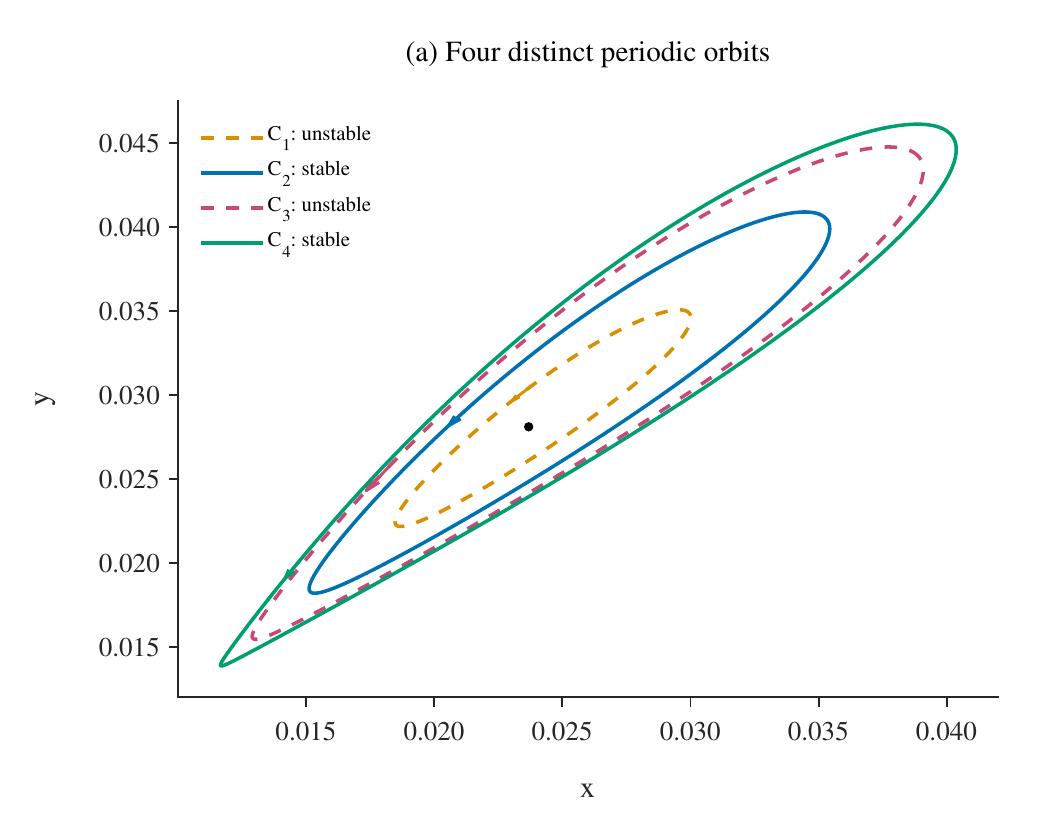}\hfill
\includegraphics[width=0.49\textwidth]{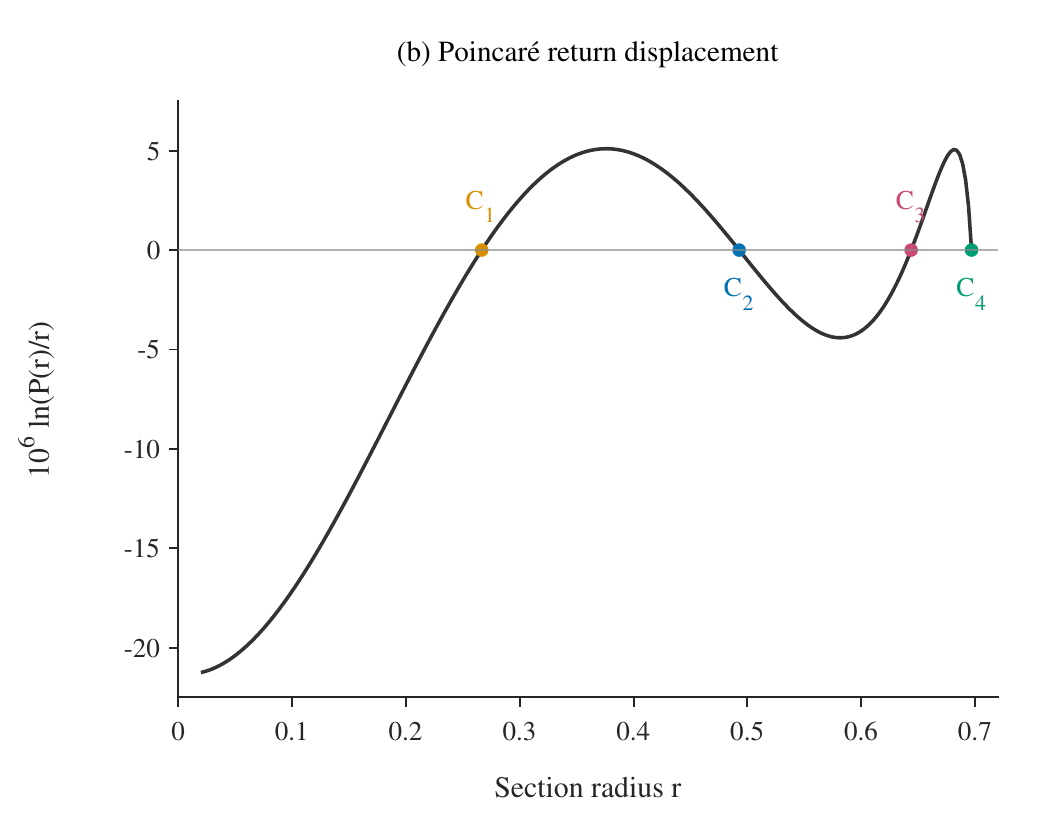}
\caption{Four coexisting limit cycles at the parameter values in \eqref{eq:numerical-parameters}, with $\beta=\beta_0$. (a) Periodic orbits in the original coordinates; the black point denotes $E_0$, dashed curves are unstable, and solid curves are stable. (b) The logarithmic return displacement $10^6L(r)$ on $\Sigma$. Its four zeros correspond to $C_1,\ldots,C_4$; the factor $10^6$ is used to make the small displacement visible.}
\label{fig:four-cycles}
\end{figure}

At $h=h_0$, the four computed limit cycles surround $E_0$ and lie in the positive quadrant; see Figure~\ref{fig:four-cycles}(a). To examine the associated Poincar\'e return map, introduce the coordinates
\begin{equation}\label{eq:normalized-coordinates}
	U=\frac{x-z_0}{z_0},\qquad
	V=\frac{\ell_0U+B_0(y/y_0-1)}{\omega_0},
\end{equation}
where
\[
A_0=1-2z_0-\frac{a_0y_0}{(a_0+z_0)^2},\qquad
B_0=-\left(1-z_0-\frac{h_0}{z_0}\right),
\]
\[
\ell_0=\frac{A_0+\delta_0}{2},
\qquad
\omega_0=\sqrt{-\ell_0^2-\delta_0B_0}>0.
\]
Let $P(r)$ denote the first-return map on
\[
\Sigma=\{(U,V):V=0,\ U>0\},
\]
parameterized by $r=U$, and define the logarithmic return displacement by
\begin{equation}\label{eq:return-displacement}
	L(r):=\ln\frac{P(r)}{r}.
\end{equation}
Numerical integration of the return-map equation gives the graph in Figure~\ref{fig:four-cycles}(b), which has four distinct positive zeros.

The nontrivial Floquet multiplier is also computed independently by integrating the divergence along each periodic orbit in the original time variable. If $(x_j(t),y_j(t))$ parameterizes the periodic orbit $C_j$ with period $T_j$, then
\begin{equation}\label{eq:floquet}
	\mu_j=\exp\left\{
	\int_0^{T_j}
	\left[
	1-2x_j(t)-\frac{a_0y_j(t)}{(a_0+x_j(t))^2}
	+\delta_0-\frac{2\beta_0y_j(t)}{x_j(t)}
	\right]d t
	\right\}.
\end{equation}

The numerical results are summarized in Table~\ref{tab:cycles}. The cycles are ordered by increasing section radius, so the stability sequence in the table is the stability sequence from the innermost to the outermost orbit. At a fixed point $r_j$ of $P$, the identity
\[
\mu_j=P'(r_j)=1+r_jL'(r_j)
\]
relates the Floquet multiplier to the slope of the return-displacement curve in Figure~\ref{fig:four-cycles}(b). Thus, $C_1$ and $C_3$ are unstable, whereas $C_2$ and $C_4$ are asymptotically stable.

\begin{table}[!ht]
\centering
\caption{Section radii, original-time periods, and nontrivial Floquet multipliers of the four computed cycles. The displayed values are rounded.}
\label{tab:cycles}
\begin{tabular}{ccccc}
\toprule
Cycle & $r_j$ & $T_j$ & $\mu_j$ & Stability\\
\midrule
$C_1$ & $0.266669123004$ & $26.3200911308$ & $1.000023284719$ & Unstable\\
$C_2$ & $0.492740289327$ & $28.1972028423$ & $0.999964893450$ & Stable\\
$C_3$ & $0.643796213321$ & $32.1353726916$ & $1.000094160670$ & Unstable\\
$C_4$ & $0.696840000024$ & $40.7866259362$ & $0.999247162291$ & Stable\\
\bottomrule
\end{tabular}
\end{table}

As an independent numerical check, we integrate the original system and its divergence equation using 40-digit working precision, with accuracy and precision goals of 28 digits. Using the unrounded section radii, the largest computed closure residual in the maximum norm is less than $2\times10^{-24}$. Moreover, the return displacements at $r_j\pm10^{-6}$ have opposite signs for each $j$, in agreement with the stability classification in Table~\ref{tab:cycles}. Since all four multipliers are close to one, the transverse attraction or repulsion is weak; the closed curves in the phase portrait should therefore not be interpreted as evidence of rapid convergence. These computations provide high-precision numerical evidence for the four-cycle configuration, but they do not constitute an interval-arithmetic proof of existence.

\section{Conclusions}\label{sec:conclusion}

We have resolved the weak-focus-order question left open in the analysis of the constant-yield harvested Holling--Tanner model by Xiang, Lu, and Huang \cite{xiang2022degenerate}. Their work established a codimension-three Hopf bifurcation and the bounds $3\leq m_{\max}\leq4$, while also identifying a codimension-four Bogdanov--Takens bifurcation. We have proved that the upper bound for the Hopf weak-focus order is attained. More precisely, we construct an admissible positive Hopf equilibrium at which the first three focal values vanish and the fourth is nonzero. We further establish the full-rank condition for the three focal values on the Hopf surface and show that the trace-zero condition is transverse to the remaining parameter direction. It follows that the maximal weak-focus order is exactly four and that the original four-parameter family admits a codimension-four Hopf bifurcation. This is a distinct result from the codimension-four Bogdanov--Takens bifurcation in \cite{xiang2022degenerate}, since the two degeneracies have different spectral character and require different nondegeneracy conditions.

The proof provides an exact algebraic realization of the fourth-order degeneracy. Starting from the elimination decomposition of the focal-value variety, we isolate a real zero of the univariate polynomial $R_3$ in a rational interval. A polynomial remainder sequence then reconstructs the prey-equilibrium coordinate, and a second remainder sequence recovers the harvesting parameter. Exact reductions modulo $R_3$, together with resultant and greatest-common-divisor tests, verify the focal-value identities and the nonvanishing conditions required in the reconstruction. Rational interval estimates certify that the reconstructed point lies in the admissible parameter domain. The elimination tools themselves, including resultants and pseudo-division, were already used in \cite{xiang2022degenerate}; the contribution here is their use in a complete exact certification chain that establishes attainability, admissibility, and transversality for a fourth-order Hopf point in the original parameter family. This strategy may be useful for other polynomial or rational planar systems whose focal-value equations admit tractable elimination and whose admissibility conditions can be controlled by exact real-algebraic methods.

The local bifurcation consequence is that a suitable generic unfolding can produce four small-amplitude limit cycles near the fourth-order weak focus. Separately, our numerical continuation study detects three folds of periodic orbits as the harvesting rate varies. At a selected parameter set, the Poincar\'e return map and the nontrivial Floquet multipliers indicate four nested periodic orbits with alternating stability: the innermost and third cycles are unstable, while the second and outermost cycles are attracting. The example therefore exhibits two coexisting attracting oscillatory states in the same planar system, separated by unstable periodic orbits. It illustrates how a high-order Hopf degeneracy can be accompanied by a nontrivial organization of periodic orbits and oscillatory multistability. The cycle-fold locations and the four-cycle configuration are supported by high-precision numerical computations; they are not claimed as interval-arithmetic existence proofs.

Several questions remain open. The present analysis determines the maximal weak-focus order but does not determine the maximum number of limit cycles over the entire admissible parameter domain. It also does not classify the global bifurcations of the four computed cycles or determine the basins of attraction of the two stable cycles. A natural next step is to validate the four-cycle configuration by interval methods and to investigate the global continuation of the periodic orbits, including possible homoclinic or other large-amplitude bifurcations. More generally, applying exact algebraic certification to other parameterized planar systems may help distinguish algebraically feasible high-order degeneracies from those that are genuinely realizable under the model's admissibility constraints.

\section*{Funding}
Z. Wen was supported by the National Natural Science Foundation of China (Nos. 12071162 and 11701191), the Natural Science Foundation of Fujian Province (No. 2025J01168), Quanzhou Municipal Special Program for Key Sci-tech Talents (No. 2026QZC29GG), and the Fundamental Research Funds for the Central Universities (No. 26JCPY004). Y. Xia was supported by the National Natural Science Foundation of China under Grant (No. 12571165) and Natural Science Foundation of Guangdong Province (2026A1515011073).

\section*{Acknowledgements}
The authors used OpenAI's ChatGPT to assist with English-language editing. All ideas and proofs were developed independently by the authors.

\section*{Competing interests}
The authors declare no competing interests.

\appendix
\section{Polynomials $g_1(z,\delta)$ and $R_3(\delta)$ used in the focal-value analysis}\label{A}
\begin{align*}
&g_1(z,\delta)=61857792 z^{14}+306302976 \delta  z^{13}-209737728 z^{13}+631426560 \delta ^2 z^{12}-898392576 \delta  z^{12}+287794944 z^{12}\\
&+684649728 \delta ^3 z^{11}-1502324736 \delta ^2 z^{11}+991213824 \delta  z^{11}-186534144 z^{11}+378892928 \delta ^4 z^{10}-1108071552 \delta ^3 z^{10}\\
&+1078970816 \delta ^2 z^{10}-368699136 \delta  z^{10}+27195264 z^{10}+30296896 \delta ^5 z^9-45589696 \delta ^4 z^9-112373632 \delta ^3 z^9\\
&+305355936 \delta ^2 z^9-223799232 \delta  z^9+47656224 z^9-113540832 \delta ^6 z^8+610416288 \delta ^5 z^8-1345952176 \delta ^4 z^8\\
&+1566351952 \delta ^3 z^8-1011785696 \delta ^2 z^8+331912848 \delta  z^8-41364720 z^8-94615248 \delta ^7 z^7+574813840 \delta ^6 z^7\\
&-1463636448 \delta ^5 z^7+2019025680 \delta ^4 z^7-1623260816 \delta ^3 z^7+753383632 \delta ^2 z^7-181671888 \delta  z^7+17343072 z^7\\
&-42112152 \delta ^8 z^6+289703336 \delta ^7 z^6-854917932 \delta ^6 z^6+1403691128 \delta ^5 z^6-1390652832 \delta ^4 z^6+842481672 \delta ^3 z^6\\
&-300918540 \delta ^2 z^6+56970696 \delta  z^6-4441968 z^6-12215972 \delta ^9 z^5+93865812 \delta ^8 z^5-315710048 \delta ^7 z^5\\
&+604898186 \delta ^6 z^5-719315672 \delta ^5 z^5+542687090 \delta ^4 z^5-255583704 \delta ^3 z^5+71573366 \delta ^2 z^5-10906428 \delta  z^5\\
&+719658 z^5-2391208 \delta ^{10} z^4+20313144 \delta ^9 z^4-76764756 \delta ^8 z^4+168488021 \delta ^7 z^4-234926248 \delta ^6 z^4\\
&+213839297 \delta ^5 z^4-125983814 \delta ^4 z^4+46329823 \delta ^3 z^4-10070836 \delta ^2 z^4+1237323 \delta  z^4-71010 z^4\\
&-308932 \delta ^{11} z^3+2888531 \delta ^{10} z^3-12168173 \delta ^9 z^3+30210407 \delta ^8 z^3-48516824 \delta ^7 z^3+52057956 \delta ^6 z^3\\
&-37245723 \delta ^5 z^3+17237865 \delta ^4 z^3-4864958 \delta ^3 z^3+779391 \delta ^2 z^3-73302 \delta  z^3+3762 z^3-24699 \delta ^{12} z^2\\
&+254676 \delta ^{11} z^2-1193837 \delta ^{10} z^2+3331095 \delta ^9 z^2-6089071 \delta ^8 z^2+7571548 \delta ^7 z^2-6444147 \delta ^6 z^2\\
&+3681299 \delta ^5 z^2-1342212 \delta ^4 z^2+283130 \delta ^3 z^2-29246 \delta ^2 z^2+1536 \delta  z^2-72 z^2-1078 \delta ^{13} z+12407 \delta ^{12} z\\
&-64900 \delta ^{11} z+202518 \delta ^{10} z-417048 \delta ^9 z+592800 \delta ^8 z-590924 \delta ^7 z+410710 \delta ^6 z-193298 \delta ^5 z+57753 \delta ^4 z\\
&-9552 \delta ^3 z+612 \delta ^2 z-19 \delta ^{14}+252 \delta ^{13}-1480 \delta ^{12}+5132 \delta ^{11}-11755 \delta ^{10}+18776 \delta ^9-21448 \delta ^8+17624 \delta ^7\\
&-10297 \delta ^6+4140 \delta ^5-1072 \delta ^4+156 \delta ^3-9 \delta ^2,\\
&R_3(\delta)={}18900179842477099144704\delta^{34}
-65488540568200541003120128\delta^{33}\\
&+3538284200965892470336701312\delta^{32}
-64937675397561940036410336768\delta^{31}\\
&+426779614729368520547919051520\delta^{30}
+212467518954161699827836364544\delta^{29}\\
&-12972716360151266948830290330432\delta^{28}
+52345886169459438466025341112304\delta^{27}\\
&-64654424068966938664601309932328\delta^{26}
-7649282526940488350149191999896\delta^{25}\\
&-533697780109492757958003790388260\delta^{24}
+3760827149808317811609859492426494\delta^{23}\\
&-9959186568305541040091171190248684\delta^{22}
+11207203014587091938650920626964801\delta^{21}\\
&+5136072942381009484139070413609345\delta^{20}
-35126480884478921220998829802237803\delta^{19}\\
&+47195309130314992941008117866806245\delta^{18}
-15622823182918656669821060643122219\delta^{17}\\
&-35793536613131613062938894770735369\delta^{16}
+54044956619105658131792376225095655\delta^{15}\\
&-25292401364755101549344397142813749\delta^{14}
-10441306745728835151971455531845111\delta^{13}\\
&+18392836043882150628155786139945939\delta^{12}
-7113938222499534650273338888719921\delta^{11}\\
&-215768623231514121675235836022917\delta^{10}
-774966293982851324432537724257421\delta^9\\
&+2187072087780173451551468544045069\delta^8
-1127668014018555291439935504570117\delta^7\\
&-103172715493790840901386443424943\delta^6
+334050104883845951143983608384790\delta^5\\
&-150177146637669865997408097581700\delta^4
+32784862916492193294882878701500\delta^3\\
&-3865060316851870806460405575000\delta^2
+235849525723928899706400000000\delta\\
&-5847857065635858368064000000.
\end{align*}

%\bibliographystyle{elsarticle-num-nodoi}
%\bibliography{RefPredator}

\end{document}